\documentclass[a4paper, 11pt]{article}
\usepackage{amsmath}
\usepackage{amsfonts}
\usepackage{amssymb}
\usepackage[english]{babel}
\usepackage{graphicx}
\usepackage{amsthm}
\usepackage[title]{appendix}

\usepackage{longtable}
\usepackage{tabularx}
\allowdisplaybreaks

\newcommand{\C}{{\cal C}}

\newcommand{\Hh}{{\cal H}}

\newtheorem{theorem}{Theorem}[section]

\newtheorem{lemma}[theorem]{Lemma}

\def\whitebox{{\hbox{\hskip 1pt
 \vrule height 6pt depth 1.5pt
 \lower 1.5pt\vbox to 7.5pt{\hrule width
    3.2pt\vfill\hrule width 3.2pt}%
 \vrule height 6pt depth 1.5pt
 \hskip 1pt } }}
\def\qed{\ifhmode\allowbreak\else\nobreak\fi\hfill\quad\nobreak
     \whitebox\medbreak}

\newcommand{\ignore}[1]{}
\usepackage[numbers,sort&compress]{natbib}

\begin{document}
\baselineskip 16pt

\title{ Almost resolvable $k$-cycle systems with $k\equiv 2\pmod 4$ }

\author{\small  L. Wang$^{1,2}$ and H. Cao$^1$ \thanks{Research
supported by the National Natural Science Foundation of China
under Grant 11571179, and the Priority Academic
Program Development of Jiangsu Higher Education Institutions.
E-mail: {\sf caohaitao@njnu.edu.cn}} \\
\small $^{1}$ Institute of Mathematics, Nanjing Normal University, Nanjing 210023, China \\
\small $^{2}$ Department of Literature and Science, Suqian College, Suqian 223800, China}

\date{}
\maketitle
\begin{abstract}

In this paper, we show that if $k\geq 6$ and $k \equiv 2 \pmod 4$, then there exists an almost resolvable $k$-cycle system of order $2kt+1$ for all $t\ge 1$
except possibly for $t=2$ and $k\geq 14$. Thus we give a partial solution to an open problem posed by Lindner, Meszka, and Rosa
(J. Combin. Des., vol. 17, pp.404-410, 2009).

\smallskip
\end{abstract}

\section{Introduction}

In this paper, we use $V(H)$ and $E(H)$ to denote the vertex-set and
the edge-set of a graph $H$, respectively.   We  denote the
cycle of length $k$ by $C_k$ and the complete graph on $v$ vertices
by $K_v$.
A {\it factor} of a graph $H$ is a spanning
subgraph whose vertex-set coincides with $V(H)$.
If its connected components are isomorphic to $G$, we call it a {\it
$G$-factor}.
A {\it $G$-factorization} of $H$ is a set of edge-disjoint $G$-factors of $H$
whose edge-sets partition $E(H)$.
An $r$-regular factor is called an {\it $r$-factor}.
Also, a {\it $2$-factorization} of a graph $H$ is a partition of $E(H)$
into 2-factors.

A $k$-{\it cycle system} of order $v$ is  a collection of
$k$-cycles whose edges partition $E(K_v)$.
A $k$-cycle system of order $v$ exists if and only if $ 3 \leq k \leq v$, $v \equiv 1
\pmod 2$ and $v(v-1) \equiv 0 \pmod {2k}$ \cite{AG,B-2003,SM,WB}.
A $k$-cycle system of order $v$ is {\it resolvable} if it has a $C_k$-factorization.
A resolvable $k$-cycle system of order $v$
exists if and only if $3 \leq k \leq v$, $v$ and $k$ are odd, and
$v \equiv 0 \pmod k$, see \cite{AH, ASSW}. If $v \equiv 1 \pmod {2k}$, then a
$k$-cycle system exists, but it is not resolvable. In this case, Vanstone et al.
\cite{VSS} started the research of the existence of almost
resolvable $k$-cycle systems.

In a $k$-cycle system of order $v$, a collection of $(v-1)/k$ disjoint $k$-cycles is called an {\it
almost parallel class}. In a $k$-cycle system of order $v \equiv 1
\pmod {2k}$, the maximum possible number of almost parallel
classes is $(v-1)/2$, in which case a  half-parallel class
containing $(v-1)/2k$ disjoint $k$-cycles is left over.
A $k$-cycle system of order $v$ whose cycle set can be partitioned
into $(v-1)/2$ almost parallel classes and a half-parallel class
is called an {\it almost resolvable $k$-cycle system}, denoted by
$k$-ARCS$(v)$.  Lindner, Meszka, and Rosa \cite{LMR} (also, Adams et al.  \cite{ABHL}) presented the following open problem
``The outstanding problem remains the construction of almost resolvable $2k$-cycle
systems of order $4k+1$, since this will determine the spectrum for almost resolvable $2k$-cycle
systems with the one possible exception of orders $8k+1$".
Since then, many authors have contributed to proving the following known conclusions.

\begin{theorem} {\rm (\cite{ABHL,BHL,CNT,DLM,DLR,LMR,NC,VSS,WC-2,WLC})}
\label{3-14} Let $n\equiv 1\pmod{2k}$. There exists a $k$-{\rm ARCS}$(n)$ for any odd $k\ge 3$
and any even $k \in \{ 4,  6,
8,  10, 14\}$, except for $(k, n) \in \{(3, 7), (3, 13), (4,
9)\}$ and except possibly for $(k, n) \in \{(8, 33), (14, 57)\}$.
\end{theorem}

In this paper, we construct
almost resolvable cycle systems of order $2k+1$ for any $k\ge 18$ and $k \equiv 2 \pmod 4$ by using algebraic methods.
Thus we have partially solved the above open problem  given by Lindner et al. in \cite{ABHL,LMR}.
Combining the known results in Theorem~\ref{3-14}, we
will prove the following main result.

\begin{theorem}
\label{main-2k+1} If $ k\geq 6$ and $k \equiv 2 \pmod 4$, then there exists a $k$-{\rm ARCS}$(2kt+1)$ for all $t\geq 1$
except possibly for $t=2$ and $k\geq 14$.
\end{theorem}

\section{Preliminary}

In this section, we present some preliminary notation and definitions, and provide lemma for the construction
of a $k$-ARCS$(2k+1)$ for $k \equiv 2 \pmod 4$.
We point out that similar methods have been used for many years (see \cite{BR,R,T2013,WC-2,WLC}).

Suppose $\Gamma$ is an additive group and  $I=\{\infty_1,\infty_2,\ldots,\infty_f\}$ is a set
which is disjoint with $\Gamma$.
We will consider an action of
$\Gamma$ on $\Gamma \ \cup \ I$ which coincides with the
right regular action on the elements of $\Gamma$, and the action of $\Gamma$ on $I$ will coincide with the identity map.
In other words, for any $\gamma\in \Gamma$, we have that $x+\gamma$ is the image under $\gamma$ of any
$x \in \Gamma$, and $x+\gamma=x$ holds for any $x \in I$.
Given a graph $H$ with vertices in $\Gamma\ \cup\ I$, the \emph{translate}
of $H$ by an element $\gamma$ of $\Gamma$ is the graph $H+\gamma$
obtained from $H$ by replacing each vertex $x\in V(H)$ with the
vertex $x+\gamma$. The {\it development} of $H$ under a subgroup
$\Sigma$ of  $\Gamma$ is the collection
$dev_\Sigma(H)=\{H+x\;|\;x\in \Sigma\}$ of all translates of $H$
by an element of $\Sigma$.

For our constructions, we set $\Gamma=\mathbb{Z}_{l}\times \mathbb{Z}_{4}$.
Given a graph $H$ with vertices in $\Gamma$ and any pair $(r,s)\in \mathbb{Z}_{4}\times
\mathbb{Z}_{4}$, we set $\Delta_{(r,s)}H=\{x-y \ | \  \{(x,r),(y,s) \} \in E(H)\}$.
Finally, given a list $\Hh = \{H_1,H_2,\ldots,H_t\}$ of graphs, we denote by
$\Delta_{(r,s)} \Hh=\cup_{i=1}^t \Delta_{(r,s)} H_i$ the multiset union of the $\Delta_{(r,s)} H_i$s.

\begin{lemma}\label{2k+1-A}
Let $v=2k+1$, $k \equiv 2 \pmod 4$, and $\C =\{F_1, F_2\}$ where each $F_i\ (i=1,2)$ is a vertex-disjoint
union of two cycles of length $k$ satisfying the following conditions:\\
$(i)$ $V(F_i)=((\mathbb{Z}_{\frac{k}{2}}\times \mathbb{Z}_{4})\cup\{\infty\})\backslash \{(a_i,b_i)\}$ for
 some $(a_i,b_i)\in \mathbb{Z}_{\frac{k}{2}}\times \mathbb{Z}_{4}$, $i=1,2$;\\
$(ii)$ $\infty$  has a neighbor in $\mathbb{Z}_{\frac{k}{2}}\times\{j\}$ for each $j \in \mathbb{Z}_{4}$;\\
$(iii)$ $\Delta_{(p, p)}\C = \mathbb{Z}_{\frac{k}{2}}\setminus\{0\}$ for each $p \in \mathbb{Z}_4$;\\
$(iv)$ $\Delta_{(0,2)}\C = \Delta_{(2,0)}\C = \mathbb{Z}_{\frac{k}{2}}\setminus\{\pm d\}$,
where $d$ satisfies $(d,\frac{k}{2})=1$;\\
$(v)$ $\Delta_{(r, s)}\C = \mathbb{Z}_{\frac{k}{2}}$ for each pair $(r,s) \in \mathbb{Z}_{4}\times \mathbb{Z}_{4}$
satisfying $ r \neq s$ and $(r,s) \notin \{ (0,2), (2,0)\}$.\\
Then, there exists a $k$-{\rm ARCS}$(v)$.
\end{lemma}

\noindent{\it Proof:\ } Let $V(K_{v})=(\mathbb{Z}_{\frac{k}{2}}\times
\mathbb{Z}_{4})\cup\{\infty\}$.  Note that $0,d,2d,\ldots,(\frac{k}{2}-1)d$ are $\frac{k}{2}$ distinct elements
since $(d,\frac{k}{2})=1$.
Then we have the required half parallel class which is formed by the cycle
$((0,0),(d,2)$, $(2d,0),(3d,2), \ldots, ((\frac{k}{2}-3)d,0),((\frac{k}{2}-2)d,2), ((\frac{k}{2}-1)d,0),$
$ (0,2),(d,0)$, $(2d,2),(3d,0), \ldots, ((\frac{k}{2}-3)d,2),((\frac{k}{2}-2)d,0), ((\frac{k}{2}-1)d,2))$.
By $(i)$, we know that $F_i$ is an almost parallel class. All the required $k$ almost parallel
classes are $F_i+(l,0)$, $i=1,2$, $l \in \mathbb{Z}_{\frac{k}{2}}$.

Now we show that the half parallel class and the $k$ almost
parallel classes form a $k$-ARCS$(v)$.
 Let  $F'$ be a graph with the edge-set $\{\{(a,0),(a+d,2)\}, \{(a,0),(a-d,2)\} \ | \ a\in \mathbb{Z}_{\frac{k}{2}} \}$ and
 $\Sigma:=\mathbb{Z}_{\frac{k}{2}}\times\{0\}$. Let $\mathcal{F} = dev_{\Sigma}(\C)\ \cup \ F'$.
The total number of edges -- counted with their respective
multiplicities -- covered by the almost parallel classes and half parallel class of $\cal F$ is
$k(2k+1)$, that is exactly the size of $E(K_v)$. Therefore, we only
need to prove that every pair of vertices  lies in a suitable
translate of $\C$ or in $F'$. By $(ii)$,  an edge $\{(z, j),
\infty\}$ of ${K}_v$ must appear in a cycle of $ dev_{\Sigma}(\C)$.

Now consider an edge $\{(z, j), (z', j')\}$ of ${K}_v$ whose
vertices both belong to $\mathbb{Z}_{\frac{k}{2}}\times \mathbb{Z}_{4}$. If $(j,j')\in \{(0,2), (2,0) \}$ and $z-z'
\in \{ \pm d\}$, then this edge belongs to $F'$. In all other cases there
is,  by $(iii)$-$(v)$, an edge of some $F_i$ of the form $\{(w, j), (w', j')\}$
such that $w-w'=z-z'$. It then follows that $F_i+(-w'+z', 0)$ is an
almost parallel class of $ dev_{\Sigma}(F_i)$ containing the edge
$\{(z, j), (z', j')\}$ and the proof is complete. \qed

\section{ Main result}

We first explain a notion which will be used in the proof of our construction. If a cycle $C$ is the concatenation of the paths
$T_0, T_1,\dots, T_m$ each of which can be obtained from a general formula, then we define $C=(T_0, \dots,\underline{T_i},\dots,T_m)$, $0\le i\le m$.  For example, $T=((-2,0),(2,3),\ldots,\underline{(-(2+i),0),(2+i,3)},\ldots,(-n,0),(n,3))$, $0\leq i \leq n-2$, means that $T$ can be viewed as  the concatenation of the paths
$T_0, T_1,\dots, T_{n-2}$, where the general formula is $T_i=((-(2+i),0),(2+i,3))$, $0 \leq i \leq n-2$.

\begin{lemma}
\label{4k+1} For any $ k\geq 18$ and $k \equiv 2 \pmod 4$, there exists a $k$-{\rm ARCS}$(2k+1)$.
\end{lemma}

\noindent {\it Proof:}  Let $v=2k+1$, and $k=4n+2$, $n\ge 4$.
We use Lemma~\ref{2k+1-A} to construct a
$k$-{\rm ARCS}$(v)$ with $V(K_{v})=(\mathbb{Z}_{\frac{k}{2}}\times
\mathbb{Z}_{4})\cup\{\infty\}$.
Three of the required parameters in $(i)$ and $(iv)$ of Lemma~\ref{2k+1-A} are
$(a_1,b_1)=(n,0)$ and $d=2$.
The other required parameters $a_2$, $ b_2$, and four cycles $\{C_1, C_2\}$ in $F_1$ and $\{C_3, C_4\}$ in $F_2$ for each $n$
are listed as below.

The cycle $C_1$ is the concatenation of the paths
$T_1$, $T_2$, $T_3$, and $T_4$, where\\
{\footnotesize
$T_1=((0,3),(0,0),(1,3))$;\\
$T_2=((-2,0),(2,3),\ldots,\underline{(-(2+i),0),(2+i,3)},\ldots,(-n,0),(n,3))$, $0\leq i \leq n-2$;\\
$T_3=((-n,1),(n,2),\ldots,\underline{(-(n-i),1),(n-i,2)},\ldots,(-2,1),(2,2))$, $0\leq i \leq n-2$;\\
$T_4=((-1,1),(0,2),(0,1))$.
}

The cycle $C_2$ is the concatenation of the paths $\infty$,
$T_1$, $T_2$, $T_3$, and $T_4$, where\\
{\footnotesize
$T_1=((1,1),(-1,3),\ldots,\underline{(1+i,1),(-(1+i),3)},\ldots,(n,1),(-n,3))$, $0\leq i \leq n-1$;\\
$T_2=((n-1,0),(-1,0),(-2,2))$;\\
$T_3=((1,0),(-3,2),\ldots,\underline{(1+i,0),(-(3+i),2)},\ldots,(n-2,0),(-n,2))$, $0\leq i \leq n-3$;\\
$T_4=((1,2),(-1,2))$.
}

\vspace{5pt}

To construct $C_3$ and $C_4$, we start with $k=18$. Here, $(a_2,b_2)=(3,2)$.\\
{\footnotesize
$C_3=(( 0, 0),( 0, 1),( 1, 0),( 2, 1),( 4, 0),( 1, 1),(-4, 0),(-2, 1),( 3, 0),( 4, 2),( 0, 3),( 1, 2),( 1, 3),(-3, 2),(-1, 3),$ \\
\hspace*{0.6cm}  $( 2, 2),( 3, 3),( 0, 2))$;\\
$C_4=(\infty,(-2, 0),(-3, 0),(-1, 0),( 2, 0),(-4, 1),( 4, 1),(-1, 1),(-3, 1),( 3, 1),(-4, 2),(-1, 2),(-2, 2),(-4, 3),$ \\
\hspace*{0.6cm}  $( 2, 3),(-2, 3),(-3, 3),( 4, 3))$.
}

\vspace{5pt}

For $k>18$, the cycle $C_3$ is the concatenation of the paths
$T_1$, $T_2$, $T_3$, and $T_4$, where\\
{\footnotesize
$T_1=((0,0),(0,1),(1,0))$;\\
$T_2=((-2,1),(2,0),\ldots,\underline{(-(2+i),1),(2+i,0)},\ldots,(-n,1),(n,0))$, $0\leq i \leq n-2$;\\
$T_3=((-n,2),(n,3),\ldots,\underline{(-(n-i),2),(n-i,3)},\ldots,(-2,2),(2,3))$, $0\leq i \leq n-2$;\\
$T_4=((-1,2),(0,3),(0,2))$.
}

Next, we consider the last cycle $C_4$. We distinguish the following two cases.

{\bf Case 1:} $k\equiv 2\pmod{8}$ and $ k\geq 26$. Here, $(a_2,b_2)=(\frac{n+2}{2},2)$.

The cycle $C_4$ is the concatenation of the paths
$\infty$, $T_1$, $T_2$, $\ldots$, $T_8$, where \\
{\footnotesize
$T_1=((-\frac{n}{2},0),(-\frac{n+2}{2},0),\ldots,\underline{(-(\frac{n}{2}-i),0),(-(\frac{n+2}{2}+i),0)},\ldots,(-1,0),(-n,0))$, $0\leq i \leq \frac{n-2}{2}$.
}

For the paths $T_2$, $T_3$, $\ldots$, $T_8$, we distinguish the following three subcases.

{\bf Case 1.1:} $k\equiv 2 \pmod{24}$ and $ k\geq 26$. \\
{\footnotesize
$T_2=((n-1,1),(2,1),(n,1),(-1,1))$;\\
$T_3=((n-2,1),(3,1),(n-4,1),(5,1),(n-3,1),(1,1),\ldots, \underline{(n-2-3i,1),(3+3i,1),}\underline{(n-4-3i,1),(5+3i,1),}$\\
\hspace*{0.6cm}  $\underline{(n-3-3i,1),(1+3i,1)},\ldots,(\frac{n+8}{2},1),(\frac{n-6}{2},1),$ $(\frac{n+4}{2},1),(\frac{n-2}{2},1),(\frac{n+6}{2},1),(\frac{n-10}{2},1))$,
$0\leq i \leq \frac{n-12}{6}$;\\
$T_4=((\frac{n+2}{2},1),(\frac{n}{2},1),(\frac{n-4}{2},1),(\frac{n}{2},2))$;\\
$T_5=((\frac{n-2}{2},2),(\frac{n+4}{2},2),\ldots,\underline{(\frac{n-2}{2}-i,2),(\frac{n+4}{2}+i,2)},\ldots,(1,2),(n,2))$, $0\leq i \leq \frac{n-4}{2}$;\\
$T_6=((-(n-1),3),(-2,3),(-n,3),(1,3))$;\\
$T_7=((-(n-2),3),(-3,3),(-(n-4),3),(-5,3),(-(n-3),3),(-1,3),\ldots, \underline{(-(n-2-3i),3),(-(3+3i),3),}$ \\
\hspace*{0.6cm}  $\underline{(-(n-4-3i),3),(-(5+3i),3),}\underline{(-(n-3-3i),3),(-(1+3i),3)},\ldots,(-\frac{n+8}{2},3),(-\frac{n-6}{2},3),$ $(-\frac{n+4}{2},3),$ \\
\hspace*{0.6cm}  $(-\frac{n-2}{2},3),(-\frac{n+6}{2},3),(-\frac{n-10}{2},3))$,
$0\leq i \leq \frac{n-12}{6}$;\\
$T_8=((-\frac{n+2}{2},3),(-\frac{n}{2},3),(-\frac{n-4}{2},3))$.
}

\vspace{5pt}

{\bf Case 1.2:} $k\equiv 10 \pmod{24}$ and $ k\geq 34$. \\
{\footnotesize
$T_2=((n-1,1),(1,1),(n-2,1),(-1,1))$;\\
$T_3=((n,1),(4,1),(n-5,1),(3,1),(n-4,1),(2,1),\ldots, \underline{(n-3i,1),(4+3i,1),}\underline{(n-5-3i,1),(3+3i,1),}$\\
\hspace*{0.6cm}  $\underline{(n-4-3i,1),(2+3i,1)},\ldots,(\frac{n+14}{2},1),(\frac{n-6}{2},1),$ $(\frac{n+4}{2},1),(\frac{n-8}{2},1),(\frac{n+6}{2},1),(\frac{n-10}{2},1))$,
$0\leq i \leq \frac{n-14}{6}$;\\
$T_4=((\frac{n+8}{2},1),(\frac{n}{2},1),(\frac{n-2}{2},1),(\frac{n+2}{2},1),(\frac{n-4}{2},1),(\frac{n}{2},2))$;\\
$T_5=((\frac{n-2}{2},2),(\frac{n+4}{2},2),\ldots,\underline{(\frac{n-2}{2}-i,2),(\frac{n+4}{2}+i,2)},\ldots,(1,2),(n,2))$, $0\leq i \leq \frac{n-4}{2}$;\\
$T_6=((-(n-1),3),(-1,3),(-(n-2),3),(1,3))$;\\
$T_7=((-n,3),(-4,3),(-(n-5),3),(-3,3),(-(n-4),3),(-2,3),\ldots, \underline{(-(n-3i),3),(-(4+3i),3),}$ \\
\hspace*{0.6cm}  $\underline{(-(n-5-3i),3),(-(3+3i),3),}\underline{(-(n-4-3i),3),(-(2+3i),3)},\ldots,(-\frac{n+14}{2},3),(-\frac{n-6}{2},3),$ $(-\frac{n+4}{2},3),$ \\
\hspace*{0.6cm}  $(-\frac{n-8}{2},3),(-\frac{n+6}{2},3),(-\frac{n-10}{2},3))$,
$0\leq i \leq \frac{n-14}{6}$;\\
$T_8=((-\frac{n+8}{2},3),(-\frac{n}{2},3),(-\frac{n-2}{2},3),(-\frac{n+2}{2},3),(-\frac{n-4}{2},3))$.
}

\vspace{5pt}

{\bf Case 1.3:} $k\equiv 18 \pmod{24}$ and $ k\geq 42$. \\
{\footnotesize
$T_2=((n-1,1),(1,1),(n-2,1),(-1,1))$;\\
$T_3=((n,1),(4,1),(n-5,1),(3,1),(n-4,1),(2,1),\ldots, \underline{(n-3i,1),(4+3i,1),}\underline{(n-5-3i,1),(3+3i,1),}$\\
\hspace*{0.6cm}  $\underline{(n-4-3i,1),(2+3i,1)},\ldots,(\frac{n+16}{2},1),(\frac{n-8}{2},1),$ $(\frac{n+6}{2},1),(\frac{n-10}{2},1),(\frac{n+8}{2},1),(\frac{n-12}{2},1))$,
$0\leq i \leq \frac{n-16}{6}$;\\
$T_4=((\frac{n+10}{2},1),(\frac{n-2}{2},1),(\frac{n+4}{2},1),(\frac{n-6}{2},1),(\frac{n+2}{2},1),(\frac{n}{2},1),(\frac{n-4}{2},1),(\frac{n}{2},2))$;\\
$T_5=((\frac{n-2}{2},2),(\frac{n+4}{2},2),\ldots,\underline{(\frac{n-2}{2}-i,2),(\frac{n+4}{2}+i,2)},\ldots,(1,2),(n,2))$, $0\leq i \leq \frac{n-4}{2}$;\\
$T_6=((-(n-1),3),(-1,3),(-(n-2),3),(1,3))$;\\
$T_7=((-n,3),(-4,3),(-(n-5),3),(-3,3),(-(n-4),3),(-2,3),\ldots, \underline{(-(n-3i),3),(-(4+3i),3),}$ \\
\hspace*{0.6cm}  $\underline{(-(n-5-3i),3),(-(3+3i),3),}\underline{(-(n-4-3i),3),(-(2+3i),3)},\ldots,(-\frac{n+16}{2},3),(-\frac{n-8}{2},3),$ $(-\frac{n+6}{2},3),$ \\
\hspace*{0.6cm}  $(-\frac{n-10}{2},3),(-\frac{n+8}{2},3),(-\frac{n-12}{2},3))$,
$0\leq i \leq \frac{n-16}{6}$;\\
$T_8=((-\frac{n+10}{2},3),(-\frac{n-2}{2},3),(-\frac{n+4}{2},3),(-\frac{n-6}{2},3),(-\frac{n+2}{2},3),(-\frac{n}{2},3),(-\frac{n-4}{2},3))$.
}

{\bf Case 2:} $k\equiv 6\pmod{8}$ and $ k\geq 22$. Here, $(a_2,b_2)=(\frac{n-1}{2},2)$

The cycle $C_4$ is the concatenation of the paths
$\infty$, $(-\frac{n+1}{2},0)$, $T_1$, $T_2$, $\ldots$, $T_8$, where \\
{\footnotesize
$T_1=((-\frac{n-1}{2},0),(-\frac{n+3}{2},0),\ldots,\underline{(-(\frac{n-1}{2}-i),0),(-(\frac{n+3}{2}+i),0)},\ldots,(-1,0),(-n,0))$, $0\leq i \leq \frac{n-3}{2}$;\\
$T_2=((n-1,1),(-1,1))$.
}

For the paths $T_3$, $T_4$, $\ldots$, $T_8$, we distinguish the following three subcases.

\vspace{5pt}

{\bf Case 2.1:} $k\equiv 6\pmod{24}$ and $ k\geq 30$. \\
{\footnotesize
$T_3=((n-2,1),(2,1),(n,1),(3,1),(n-4,1),(1,1),\ldots, \underline{(n-2-3i,1),(2+3i,1),}\underline{(n-3i,1),(3+3i,1),}$\\
\hspace*{0.6cm}  $\underline{(n-4-3i,1),(1+3i,1)},\ldots,(\frac{n+9}{2},1),(\frac{n-9}{2},1),$ $(\frac{n+13}{2},1),(\frac{n-7}{2},1),(\frac{n+5}{2},1),(\frac{n-11}{2},1))$,
$0\leq i \leq \frac{n-13}{6}$;\\
$T_4=((\frac{n+3}{2},1),(\frac{n+1}{2},1),(\frac{n-5}{2},1),(\frac{n-1}{2},1),(\frac{n+7}{2},1),(\frac{n-3}{2},1),(\frac{n+1}{2},2),(\frac{n+3}{2},2))$;\\
$T_5=((\frac{n-3}{2},2),(\frac{n+5}{2},2),\ldots,\underline{(\frac{n-3}{2}-i,2),(\frac{n+5}{2}+i,2)},\ldots,(1,2),(n,2))$, $0\leq i \leq \frac{n-5}{2}$;\\
$T_6=((-(n-1),3),(1,3))$;\\
$T_7=((-(n-2),3),(-2,3),(-n,3),(-3,3),(-(n-4),3),(-1,3),\ldots, \underline{(-(n-2-3i),3),(-(2+3i),3),}$ \\
\hspace*{0.6cm}  $\underline{(-(n-3i),3),(-(3+3i),3),}\underline{(-(n-4-3i),3),(-(1+3i),3)},\ldots,(-\frac{n+9}{2},3),(-\frac{n-9}{2},3),$ $(-\frac{n+13}{2},3),$ \\
\hspace*{0.6cm}  $(-\frac{n-7}{2},3),(-\frac{n+5}{2},3),(-\frac{n-11}{2},3))$,
$0\leq i \leq \frac{n-13}{6}$;\\
$T_8=((-\frac{n+3}{2},3),(-\frac{n+1}{2},3),(-\frac{n-5}{2},3),(-\frac{n-1}{2},3),(-\frac{n+7}{2},3),(-\frac{n-3}{2},3))$.
}

\vspace{5pt}

{\bf Case 2.2:} $k\equiv 14\pmod{24}$ and $ k\geq 38$. \\
{\footnotesize
$T_3=((n-2,1),(2,1),(n,1),(3,1),(n-4,1),(1,1),\ldots, \underline{(n-2-3i,1),(2+3i,1),}\underline{(n-3i,1),(3+3i,1),}$\\
\hspace*{0.6cm}  $\underline{(n-4-3i,1),(1+3i,1)},\ldots,(\frac{n+11}{2},1),(\frac{n-11}{2},1),$ $(\frac{n+15}{2},1),(\frac{n-9}{2},1),(\frac{n+7}{2},1),(\frac{n-13}{2},1))$,
$0\leq i \leq \frac{n-15}{6}$;\\
$T_4=((\frac{n+5}{2},1),(\frac{n-7}{2},1),(\frac{n-1}{2},1),(\frac{n+9}{2},1),(\frac{n-5}{2},1),(\frac{n+3}{2},1),(\frac{n+1}{2},1),(\frac{n-3}{2},1),(\frac{n+1}{2},2),(\frac{n+3}{2},2))$;\\
$T_5=((\frac{n-3}{2},2),(\frac{n+5}{2},2),\ldots,\underline{(\frac{n-3}{2}-i,2),(\frac{n+5}{2}+i,2)},\ldots,(1,2),(n,2))$, $0\leq i \leq \frac{n-5}{2}$;\\
$T_6=((-(n-1),3),(1,3))$;\\
$T_7=((-(n-2),3),(-2,3),(-n,3),(-3,3),(-(n-4),3),(-1,3),\ldots, \underline{(-(n-2-3i),3),(-(2+3i),3),}$ \\
\hspace*{0.6cm}  $\underline{(-(n-3i),3),(-(3+3i),3),}\underline{(-(n-4-3i),3),(-(1+3i),3)},\ldots,(-\frac{n+11}{2},3),(-\frac{n-11}{2},3),$ $(-\frac{n+15}{2},3),$ \\
\hspace*{0.6cm}  $(-\frac{n-9}{2},3),(-\frac{n+7}{2},3),(-\frac{n-13}{2},3))$,
$0\leq i \leq \frac{n-15}{6}$;\\
$T_8=((-\frac{n+5}{2},3),(-\frac{n-7}{2},3),(-\frac{n-1}{2},3),(-\frac{n+9}{2},3),(-\frac{n-5}{2},3),(-\frac{n+3}{2},3),(-\frac{n+1}{2},3),(-\frac{n-3}{2},3))$.
}

\vspace{5pt}

{\bf Case 2.3:} $k\equiv 22\pmod{24}$ and $ k\geq 22$.\\
{\footnotesize
$T_3=((n-2,1),(2,1),(n,1),(3,1),(n-4,1),(1,1),\ldots, \underline{(n-2-3i,1),(2+3i,1),}\underline{(n-3i,1),(3+3i,1),}$\\
\hspace*{0.6cm}  $\underline{(n-4-3i,1),(1+3i,1)},\ldots,(\frac{n+7}{2},1),(\frac{n-7}{2},1),$ $(\frac{n+11}{2},1),(\frac{n-5}{2},1),(\frac{n+3}{2},1),(\frac{n-9}{2},1))$,
$0\leq i \leq \frac{n-11}{6}$;\\
$T_4=((\frac{n+1}{2},1),(\frac{n+5}{2},1),(\frac{n-1}{2},1),(\frac{n-3}{2},1),(\frac{n+1}{2},2),(\frac{n+3}{2},2))$;\\
$T_5=((\frac{n-3}{2},2),(\frac{n+5}{2},2),\ldots,\underline{(\frac{n-3}{2}-i,2),(\frac{n+5}{2}+i,2)},\ldots,(1,2),(n,2))$, $0\leq i \leq \frac{n-5}{2}$;\\
$T_6=((-(n-1),3),(1,3))$;\\
$T_7=((-(n-2),3),(-2,3),(-n,3),(-3,3),(-(n-4),3),(-1,3),\ldots, \underline{(-(n-2-3i),3),(-(2+3i),3),}$ \\
\hspace*{0.6cm}  $\underline{(-(n-3i),3),(-(3+3i),3),}\underline{(-(n-4-3i),3),(-(1+3i),3)},\ldots,(-\frac{n+7}{2},3),(-\frac{n-7}{2},3),$ $(-\frac{n+11}{2},3),$ \\
\hspace*{0.6cm}  $(-\frac{n-5}{2},3),(-\frac{n+3}{2},3),(-\frac{n-9}{2},3))$,
$0\leq i \leq \frac{n-11}{6}$;\\
$T_8=((-\frac{n+1}{2},3),(-\frac{n+5}{2},3),(-\frac{n-1}{2},3),(-\frac{n-3}{2},3))$.
}
\qed

Now we have constructed a $k$-ARCS$(2k+1)$ for any $k\ge 18$ and $k\equiv 2\pmod 4$. That is enough to prove our main result by using any of the two known recursive constructions
which can be found in \cite{LMR} and \cite{CNT}, respectively.

In the first recursive construction from \cite{LMR}, the authors start with a commutative quasigroup of order $2t$ with holes of size 2 (see \cite{LR}), then give each vertex weight $k$ and use a $k$-ARCS$(2k+1)$ (exists by assumption) and a $C_{k}$-factorization of the complete bipartite graph $K_{k,k}$ (see \cite{P}) as input designs to get a $k$-ARCS$(2kt+1)$ for all $t\ge 3$ and even integer $k\ge 8$.
In the second recursive construction from \cite{CNT}, the authors use a $k$-ARCS$(2k+1)$ and a $k$-cycle frame of type $(2k)^t$ (all exist by assumption) to get a $k$-ARCS$(2kt+1)$. Note that Buratti et al. \cite{BCDT} have proved that there exists a $k$-cycle frame of type $(2k)^t$ for all $t\geq 3$ when $k$ is even, and $t\geq 4$ when $k$ is odd. So for any even integer $k\ge 4$, we may use the second recursive construction to obtain a $k$-ARCS$(2kt+1)$ for all $t\ge 3$ if there is a $k$-ARCS$(2k+1)$.

Actually, a commutative quasigroup of order $2t$ with holes of size 2 and a $C_{k}$-factorization of $K_{k,k}$ can lead to a $k$-cycle frame of type $(2k)^t$ when $k\ge 8$ is even and $t\ge 3$.  Thus the two recursive constructions are the same when $k$ is even and $k\ge 8$  although they use different notations. We state it in the following theorem.

\begin{theorem} {\rm (\cite{LMR,CNT})}
\label{Con} Let $k\equiv0\pmod2$ and $k\geq8$. If there exists a {\rm $k$-ARCS$(2k+1)$},
then  there exists a {\rm $k$-ARCS$(2kt+1)$} for all $t\ge 1$ except possibly for $t=2$.
\end{theorem}

At last, we prove our main result.

\noindent {\bf  Proof of Theorem~\ref{main-2k+1}:}
Combining Theorems~\ref{3-14},~\ref{Con}, and Lemma~\ref{4k+1}, we can obtain the conclusion.
\qed

\vspace{5pt}

\noindent {\bf Acknowledgments} We would like to 
thank the anonymous referees for their careful reading of the
manuscript and many constructive comments and suggestions that
greatly improved the readability of this paper.


\begin{thebibliography}{Z}
\baselineskip 11pt


\bibitem{ABHL}
P. Adams, E. J. Billington, D. G. Hoffman, and C. C. Lindner, {The
generalized almost resolvable cycle system problem},
{\it Combinatorica} {\bf 30} (2010), 617-625.


\bibitem{AG}
B. Alspach and H. Gavlas, Cycle decompositions of $K_n$ and $K_n- I$,
 {\it J. Combin. Theory Ser. B} {\bf 81} $(2001)$, 77-99.


\bibitem{AH}
B. Alspach and R. H$\ddot{a}$ggkvist, Some observations on the Oberwolfach problem,
{\it Journal of Graph Theory }{\bf 9} $(1985)$,  177-187.


\bibitem{ASSW}
B. Alspach,  P. J. Schellenberg,  D. R. Stinson, and D. Wagner, The
Oberwolfach problem and factors of uniform odd length cycles, {\it
J. Combin.  Theory Ser.  A} {\bf 52} $(1989)$,  20-43.



\bibitem{BHL}
E. J. Billington, D. G. Hoffman, C. C. Lindner, and M. Meszka, Almost resolvable minimum coverings of complete graphs with 4-cycles,
 {\it  Australas. J. Combin.} {\bf 50} (2011), 73-85.


\bibitem{B-2003}
M. Buratti, Rotational $k$-cycle systems of order $v < 3k$; another proof of the existence of odd cycle systems,
 {\it J. Combin. Des.}  {\bf 11} $(2003)$, 433-441.


\bibitem{BCDT}
M. Buratti, H. Cao, D. Dai, and T. Traetta, A complete solution to the existence of $(k,\lambda)$-cycle frames of type $g^u$,
 {\it J. Combin. Des.} {\bf 25} $(2017)$, 197-230.


\bibitem{BR}
M. Buratti and G. Rinaldi, On sharply vertex transitive $2$-factorizations of the complete graph,
 {\it J. Combin. Theory Ser. A}  {\bf 111} $(2005)$, 245-256.


\bibitem{CNT}
H. Cao,  M. Niu, and C. Tang, On the existence of cycle frames and
almost resolvable cycle systems, {\it Discrete Math.} {\bf 311} (2011), 2220-2232.


\bibitem{DLM}
I. J. Dejter, C. C. Lindner, M. Meszka, and C. A. Rodger, Corrigendum/addendum to: almost resolvable 4-cycle systems,
 {\it J. Combin. Math. Combin. Comput. } {\bf 66} $(2008)$, 297-298.


\bibitem{DLR}
I. J. Dejter, C. C. Lindner, C. A. Rodger, and M. Meszka, Almost resolvable 4-cycle systems,
 {\it J. Combin. Math. Combin. Comput. } {\bf 63} $(2007)$, 173-181.





\bibitem{LMR}
C. C. Lindner, M. Meszka, and A. Rosa, Almost resolvable cycle systems---an analogue of Hanani triple systems,
 {\it J. Combin. Des.}  {\bf 17} $(2009)$, 404-410.



\bibitem{LR}
C. C. Lindner and C. A. Rodger, Design Theory,
CRC Press, Boca Raton, FL, 1997, 208pp.








\bibitem{NC}
M. Niu and H. Cao, More results on cycle frames and almost resolvable cycle systems, {\it Discrete Math.} {\bf
312} (2012), 3392-3405.



\bibitem{P}
W. L. Piotrowski, The solution of the bipartite analogue of the Oberwolfach problem,
{\it Discrete Math.}  {\bf 97} (1991), 339-356.





\bibitem{R}
 A. Rosa, On certain valuations of the vertices of a graph, in: Theory of Graphs, Internat. Sympos., Rome, 1966, Gordon and
Breach/Dunod, New York/Paris, 1967, pp. 349-355.


\bibitem{SM}
M. $\check{ \rm S}$ajna, Cycle decompositions: complete graphs and fixed length cycles,
{\it J. Combin. Des.}  {\bf 10} (2002), 27-78.


\bibitem{T2013}
T. Traetta, A complete solution to the two-table Oberwolfach problems,
 {\it J. Combin. Theory Ser. A} {\bf 120}  $(2013)$, 984-997.


\bibitem{VSS}
S. A. Vanstone, D. R. Stinson, P. J. Schellenberg, A. Rosa, R. Rees, C. J. Colbourn, M. W. Carter, and J. E. Carter,
Hanani triple systems, {\it Israel J. Math.}  {\bf 83} (1993), 305-319.



\bibitem{WC-2}
L. Wang and H. Cao, Completing the spectrum of almost resolvable cycle systems with odd cycle length,
{\it Discrete Math.} {\bf 341} (2018), 1479-1491.



\bibitem{WLC}
L. Wang, S. Lu, and H. Cao, Further results on almost resolvable cycle systems and the Hamilton-Waterloo problem,
{\it J. Combin. Des.}  {\bf 26} (2018), 27-47.



\bibitem{WB}
S. Wu and M. Buratti, A complete solution to the existence problem for 1-rotational $k$-cycle systems of $K_v$,
{\it J. Combin. Des.}  {\bf 17} (2009), 283-293.

\end{thebibliography}
\end{document}